\input amstex
\documentstyle{amsppt}
\magnification=1200
\pagewidth{6.5 true in}
\pageheight{9.0 true in}
\nologo
\tolerance=10000
\overfullrule=0pt

\define\opn{\operatorname}
\define\1{^{-1}}
\define\8{\infty}
\define\A{\frak A}

\define\brq{^{[q]}}
\define\bx{\bold x}

\define\hgt{\opn{ht}}

\define\ic#1{\overline{#1}}

\define\inc{\subseteq}

\redefine\l{\ell}

\define\ra{\longrightarrow}
\define\m{{\bold m}}

\redefine\vec#1 #2{#1_1,\ldots,#1_{#2}}

\topmatter
\title
F-rational rings and the integral closures of ideals
\endtitle
\date October 2, 2000\enddate
\rightheadtext{F-rational rings}
\author
Ian M. Aberbach and
Craig Huneke
\endauthor
\address
Department of Mathematics, University of Missouri,
Columbia, MO 65211
\endaddress
\email
aberbach\@math.missouri.edu
\endemail
\address
Department of Mathematics, University of Kansas,
Lawrence, KS 66045
\endaddress
\email
huneke\@math.ukans.edu
\endemail
\thanks
Both authors were partially supported by the NSF.
\endthanks
\endtopmatter
\document
\head 1. Introduction \endhead
\bigskip

The history of the Brian\c con-Skoda theorem and its ensuing avatars
in commutative algebra has been well-documented in many papers. For example
see \cite{LS, AH1}. We will therefore only briefly review the relevant
concepts and theorems. First recall the definitions of the integral closure of
an ideal.

\definition{Definition 1.1}
Let $R$ be a ring and let $I$ be an ideal of $R$.  An element $x\in R$ is
{\it integral over $I$\/} if $x$ satisfies an equation of the form $x^n+a_1x^{
n-1}+\dots+a_n=0$, where $a_j\in I^j$ for $1\leq j\leq n$.
The \it integral closure \rm of $I$, denoted by $\overline I$,
is the set of all elements integral over $I$. This set is an ideal.
\enddefinition

Let $R^o$ be the set of all elements of $R$ not in a minimal prime.
An equivalent though less standard (but for our purposes a more useful)
 definition of integral closure is:

\definition{Equivalent Definition 1.1} Let $R$ be a Noetherian ring and let $I$ be an ideal of $R$.  An element $x\in R$ is
{\it integral over $I$\/} if there exists an element $c\in R^o$ such that
$cx^n\in I^n$ for all $n >> 0$.
\enddefinition

A theorem proved  by Brian\c con and Skoda \cite{BS} for convergent power
series over the complex numbers and generalized to arbitrary regular
local rings by Lipman and Sathaye states:
 
\proclaim{Theorem 1.2 \cite{BS,LS}} Let $R$ be a regular local ring and let
$I$ be an ideal generated by $\ell$ elements. Then for all
$n\geq \ell$, $$\overline{I^n}\subseteq I^{n-\ell+1}.$$
\endproclaim

This was partially extended to the class of pseudo-rational rings by Lipman and
Teissier \cite{LT}. However, they were unable to recover the full
strength of (1.2).

\proclaim{Theorem 1.3 (\cite{LT, (2.2)})} Let $R$ be a Noetherian
local ring and assume that the localization $R_P$ is pseudo-rational
for every prime ideal $P$ in $R$. Suppose that $I$ has a reduction
$J$ such that $\dim R_P\leq \delta$ for every associated prime
$P$ of $J^n$. Then
$$
\overline{I^{n+\delta-1}}\subseteq J^n.
$$
In particular, if $J$ can be generated by a regular sequence
of length $\delta$, then the above containment holds for all $n\geq 1$.
\endproclaim

The present two authors, as well as Lipman, have pushed the original
theorem further by introducing `coefficients'. See \cite{AH1,2, AHT, L1}.
The methods used by the present authors have relied on the theory of
tight closure. These improvements, however, have been valid only in regular
rings, and the question of whether the statement of Theorem 1.2 remains
valid in arbitrary pseudo-rational rings has remained open since 1981.
Recent progress was made by  Hyry 
and Villamayor \cite{HV}, who proved (among other things) that
if $R$ is local Gorenstein and essentially of finite type over a
field of characteristic $0$, then $\overline{I^{n+\ell-1}}\inc I^n$ for
an arbitrary ideal $I$ with $\ell$ generators.
In this paper we will use tight closure methods to prove (1.2) is
valid for F-rational rings (the definition is below). In characteristic
$p$, K. Smith \cite{Sm} proved that F-rational implies pseudo-rational, but
it can be stronger in general. However for affine algebras in equicharacteristic $0$, the
concepts of rational singularity, pseudo-rational singularity and 
F-rational type all agree, due to work of Lipman and Teissier \cite{LT}
for the equivalence of rational singularity and pseudo-rational singularity,
and of Smith \cite{Sm} and N. Hara \cite{Ha} and independently
V. Mehta and V. Srinivas \cite{MS} for the equivalence of rational singularity
and F-rational type (Smith proved that rational implies F-rational type
and the other authors have just recently proved the converse).
It follows from these equivalences that in equicharacteristic $0$ we
are able to prove (1.2) for rational singularities. 

The basic idea of this paper is inspired by the proof of a cancellation
theorem \cite{Hu}. The key idea is to relate an arbitrary ideal $I$ to
a system of parameters in a manner which closely approximates the
structure of the powers of $I$. We do this via a basic construction, and
then a theorem which relates the integral closure of  powers of $I$ with
the tight closure of the system of parameters. In the next section
 we briefly discuss tight
closure, and refer to \cite{HH1, Hu2} for more references and information.

\bigskip
\head 2. Tight Closure \endhead
\bigskip

We begin with the definition.

\definition{Definition 2.1} Let $R$ be a Noetherian ring of characteristic $p > 0$. Let $I$
be an ideal of $R$. An element $x\in R$ is said to be in the tight closure of $I$ if
there exists an element $c\in R^o$ such that for all
large $q= p^e$, $cx^q\in I^{[q]}$, where $I^{[q]}$ is the ideal generated by the $q$th
powers of all elements of $I$.
\enddefinition

Every ideal in a regular ring is tightly closed. We say elements
$x_1,...,x_n$ in $R$ are \it parameters \rm if the height of the
ideal generated by them is at least $n$ (we allow them to be
the whole ring, in which case the height is said to be $\infty$).
If the ideal they generate is proper, then the Krull height
theorem says that the height is exactly $n$.

\definition{Definition 2.2} A Noetherian ring $R$ of characteristic $p > 0$
is said to be \it F-rational \rm if the ideals generated by 
parameters are tightly closed.
\enddefinition

This definition arose from the work of Fedder and Watanabe \cite{FW} because of
the apparent connection to the concept of rational singularities.

The concept of pseudo-rationality was introduced in \cite{LT}, partly as a
substitute for the notion of rational singularities in positive and mixed
characteristic, where desingularizations are not known to exist in general.

Their definition is \cite{LT, Section 2}:

\definition{Definition 2.3} Let $(R,\m)$ be a $d$-dimensional local Noetherian
ring. $R$ is said to be \it pseudo-rational \rm if it is normal, Cohen-Macaulay,
analytically unramified, and if for every proper birational map
$\pi: W\ra X = \text{Spec}(R)$ with $W$ normal and closed fiber $E = \pi^{-1}(\m)$,
the canonical map
$$H^d_{\m}(\pi_*(\Cal O_W))= H^d_{\m}(R)\ra H^d_E(\Cal O_W)$$
is injective.
\enddefinition

In \cite{LT}, it is proved that for a local ring essentially of finite
type over a field of characteristic $0$, the notions of pseudo-rational
and rational singularity agree. In \cite{Sm}, it is shown that in
positive characteristic, F-rational implies pseudo-rational. 
Smith uses this to prove that rings of finite type over a field of characteristic $0$
which are F-rational type have rational singularities. `F-rational type' essentially
means that charactersitic $p$ models of the variety are F-rational. Precisely, we
need to introduce the idea of a model:

   If $R$ is a ring which is finitely generated over a field of
characteristic $0$, say $R = k[X_1,...,X_n]/I$, then we can choose
generators for the ideal $I$ and by collecting coefficients of those
generators find a finitely generated $\bold Z$-algebra $A\inc k$ such
that if we define $R_A = A[X_1,...,X_n]/(I\cap A[X_1,...,X_n])$, 
then $R = k\otimes_AR_A$. We call the map  $ A\ra R_A$ a family of models of $R$.
We sometimes insist that the map $ A\ra R_A$ be flat, which one can
always obtain by expanding $A$ by localizing at a single element.
A typical closed fiber of $R_A$ over $A$ is a characteristic $p$ model
of $R$.  
 
\definition{Definition 2.5} Let $R$ be a finitely generated algebra over
a field of characteristic $0$. $R$ is said to have \it F-rational type \rm
if $R$ admits a family of models $A\ra R_A$ in which a Zariski dense
set of closed fibers are F-rational. (This does not depend on the choice
of models.)
\enddefinition

The theorem in \cite{Sm} says that if $X$ is a scheme of finite type
over a field of characteristic $0$, then if $X$ has F-rational type
it has only rational singularities. Recently, the converse has been
proved by N. Hara \cite{Ha}, and independently by Mehta and
Srinivas \cite{MS}.

\head 3. F-rational Rings and Tight Closure \endhead
\bigskip

In this section we first discuss a basic construction which will play a
crucial role in the paper. Given an ideal $I$ in a Noetherian local
ring $(R,\m)$, a minimal reduction $J$ of $I$, say $J = (a_1,...,a_{\l})$,
and an integer $N$, we wish to construct an ideal $\A$, generated by parameters such that
$J\equiv \frak A$ modulo $\m^N$, and such that $\A$ is closely related to $I$
and its powers. For example, one would like $I\inc \A$, but this is in
general not possible since $I$ may not be contained in any ideal
generated by parameters. We record what we need in Proposition 3.2.
We need the following lemma from \cite{AHT, (7.2)}.

\proclaim{Lemma 3.1} Let $(R,\m)$ be a
local ring with infinite residue field
 and let $I\inc R$ be an ideal of analytic spread
$\l$.  Let $J \inc I$ be a minimal reduction of $I$.  Then there
exists a ``basic'' generating set $\vec a \l$ for $J$ such that
\roster
\item if $P$ is a prime ideal containing $I$ and $\hgt P = i \leq \l$
then $(\vec a i)_P$ is a reduction of $I_P$, and
\item $\hgt((\vec a i)I^n:I^{n+1} + I) \ge i+1$ for all $n \gg 0$.
\item If $c_i\equiv a_i$ modulo $I^2$, then (1) and (2) hold with
$c_i$ replacing $a_i$.
\endroster
\endproclaim

\demo{Proof}
The first two statements are found in \cite{AHT, Lemma 7.2}. The last
statement follows from the proof of Lemma 7.2 in \cite{AHT}. The choice
of a basic generating set only depends on the images of the $a_i$ in
the associated graded ring $G = R/I \oplus I/I^2\oplus ...$. In particular
since $c_i$ and $a_i$ have the same leading forms in $G$, (3) follows. \qed
\enddemo

\proclaim{Proposition 3.2}  Let $(R,\m)$ be an equidimensional and catenary
local ring with infinite residue field
 and let $I\inc R$ be an ideal of analytic spread
$\l$.  Let $J \inc I$ be a minimal reduction of $I$.
We assume  that $\hgt I = g$, and $J = (a_1',...,a_{\l}')$, a basic generating
set for $J$ as in Lemma 1.1.
Let $N$ and $w$ be fixed integers, and suppose that for
$g+1\leq i\leq \l$ we are given finite sets of
primes $\Lambda_i = \{Q_{ji}\}$ all containing $I$ and of height $i$. Then there exist
elements $a_1,...,a_{\l}$ and $t_{g+1},...,t_{\l}$ such that the following
hold. (We set $t_i = 0$ for $i\leq g$ for convenience).
\roster
\item $a_i\equiv a_i'$ modulo $I^2$.
\item For $g+1\leq i\leq \l$, $t_i\in \m^N$.
\item $b_1,...,b_g,b_{g+1},...,b_{\l}$ are parameters,
where $b_i = a_i + t_i$.
\item If  $R/I$ is equidimensional, the images of $t_{g+1},...,t_{\l}$  in $R/I$
are parameters.
\item There is an  integer $M$ such that $t_{i+1}\in (J_i^tI^{M}:I^{M+t})$
for all $0\leq t\leq w + \l$
where $J_i = (a_1,...,a_i)$.
\item $t_{i+1}\notin \cup_j Q_{ji}$, the union being over the primes in $\Lambda_i$.
\endroster
\endproclaim

\demo{Proof}  We choose the $a_i$ and $t_i$ inductively.  We first modify $a_1',...,a_g'$
to $a_1,...,a_g$ in such a way that these elements form parameters. We can do this
with $a_i\equiv a_i'$ modulo $I^2$ for $1\leq i\leq g$.
Suppose we have chosen $a_1,...,a_i$ and $t_1,...,t_i$
so that all of the above statements are true for these elements. Fix the
minimal primes $P_1,...,P_k$ (all necessarily of height $i$) above 
$B_i = (b_1,...,b_i)$. Divide them into two sets. Let $P_1,...,P_n$ be the
ones which contain $I$, and $P_{n+1},...,P_k$ those which don't contain $I$.
We first change $a_{i+1}'$ to an element $a_{i+1}\equiv a_{i+1}'$ modulo
$J^2$ such that $a_{i+1}\notin \cup_{j=n+1}^kP_j$. This choice is
possible as the nilradical of $J$ is the same as the nilradical of $I$.
Next choose $M_{i}$ such that the height of $I + (J_{i}I^{M_{i}}:I^{M_{i} + 1})$
is least $i+1$, and choose $M$ to be the maximum of the $M_i$.
(This is possible by Lemma 3.1.)
 This choice forces all $(J_i^tI^{M}:I^{M + t }) +I$ to be
height at least $i+1$ for all $t\geq 0$. For suppose that $(J_i^tI^{M}:I^{M + t }) +I\inc Q$,
where $Q$ is a prime of height at most $i$. Since $I\inc Q$, this forces
$(J_{i}I^{M}:I^{M + 1})\nsubseteq Q$, and after localization at $Q$
$(I^{M + 1})_Q = (J_{i}I^{M})_Q$. But this forces $(I^{M + t})_Q = (J_{i}^tI^{M})_Q$ for all
integers $t$, and so $(J_i^tI^{M}:I^{M + t })\nsubseteq Q$, a contradiction.
 Using prime avoidance,
choose
$$t_{i+1}\in \cap_{t=0}^{w+\l}(J_{i}^tI^{M}:I^{M + t})\cap \m^N\cap (\cap_{j=n+1}^kP_j)$$
and
$$t_{i+1}\notin (\cup_{j=1}^nP_j) \cup (\cup_j Q_{ji}).$$

This is possible since $I$ is contained in each of the primes in the
second line, but all these primes have height $i$, while
the height of  $I + (J_{i}^tI^{M}:I^{M + t})$ is at least $i+1$.
We set $b_{i+1} = a_{i+1}+ t_{i+1}$. We claim this choice proves
(1)-(6) for these new elements.  Our choice of $a_{i+1}$ and $t_{i+1}$ make
statements (1), (2), (5) and (6) trivial. 
To prove (3) we need only to prove that $b_{i+1}\notin \cup_{j=1}^k P_j$. If $j\leq n$,
then $a_{i+1}\in I\inc P_j$ while $t_{i+1}\notin P_j$. Hence $b_{i+1}\notin P_j$.
If $j\geq n+1$, then $a_{i+1}\notin P_j$ while $t_{i+1}\in P_j$. Again $b_{i+1}\notin P_j$,
proving (3).  Statement (4) follows from (3).  Clearly the height of $(I,b_{g+1},...,b_{i+1})$
is at least that of $b_1,...,b_{i+1}$, hence at least $i+1$. But
$(I,b_{g+1},...,b_{i+1}) = (I,t_{g+1},...,t_{i+1})$. As 
 $R$ is equidimensional and catenary, it follows that the images of the $t_j$
in $R/I$ form parameters. \qed\enddemo

 \proclaim{Theorem 3.3}  Let $(R,\m)$ be an equidimensional and catenary local ring
of characteristic $p$ having infinite residue field.  Let $I$ be an ideal of analytic
spread $\l$ and positive height $g$.  Let $J$ be a minimal reduction of $I$.
Fix $w, N \ge 0$.
Choose $a_i$ and $t_i$ as in Proposition 3.2. Set $\A = B_{\l} = (b_1,...,b_g,...,b_{\l})$.
 Then 
$$\ic{I^{\l+w}} \inc (\A^{w+1})^*.$$
\endproclaim

\demo{Proof}
Our choice of elements says that $a_1, \ldots, a_g, a_{g+1} + t_{g+1}, \ldots, a_{i+1}+t_{i+1}$
is a part of a system of parameters. Fix the notation as in Proposition 3.2.
  By our choice of
the $t_j$ we have  that for all $1\leq k\leq w + \l$, $t_jI^{M+k}
\inc J_{j-1}^kI^M$. We first claim that this implies
$$t_j^nI^{M+nk}\inc J_{j-1}^{nk}I^M$$ for all $n\geq 1$. Assume this is true for a fixed $n$,
and multiply by $t_jI^k$. We obtain that 
$(t_jI^k)(t_j^nI^{M+nk})\inc (t_jI^k)J_{j-1}^{nk}I^M$. Since
 $t_jI^{M+k}\inc J_{j-1}^kI^M$, we obtain that
$$t_j^{n+1}I^{M+(n+1)k}\inc J_{j-1}^{nk}J_{j-1}^kI^M$$
as required.    Fix $c \in I^M\cap R^0$. Note that the above containment
shows that for all $n\geq 1$
$$ct_j^{n+1}I^{(n+1)k}\inc J_{j-1}^{(n+1)k}. \tag{3.4}$$

Set $B_i = (b_1,...,b_i)$.
Let $g \le i \le \l$ and $w \ge r \ge 0$.  We show by induction that
$c^{i-g} J_i^{(i+r)q} \inc  (B_i^{r+1})\brq$.
The base case is when $i = g$ and $r \le w$ is arbitrary.  In this
case $J_g^{(g+r)q} \inc (J_g^{r+1})\brq = (B_g^{r+1})\brq$. The
first equality in the above line follows at once from \cite{HH1, proof of (5.4)}.

Assume now that we are given $r$ and $i > g$,
 and the claim is true for $i' < i$ (with
$r' \le w$ arbitrary) or $i' = i$ and $r' < r \le w$.
By our choice of $c$ and of the $t_j$, 
$$\multline
c^{i-g} J_i^{(i+r)q} \inc c^{i-g}J_i\brq J_i^{(i+r-1)q}
\inc c^{i-g} [J_g\brq J_i^{(i+r-1)q} +  a_{g+1}^q J_i^{(i+r-1)q}
+ \cdots +  a_{i}^q J_i^{(i+r-1)q}] \\
= c^{i-g-1}[cJ_g\brq J_i^{(i+r-1)q} + ca_{g+1}^q J_i^{(i+r-1)q}
+ \cdots + ca_{i}^q J_i^{(i+r-1)q}].\endmultline$$
Consider a typical term in this sum, $ca_j^qJ_i^{(i+r-1)q}$, where
$g+1\leq j\leq i$. As $b_j = a_j + t_j$, we can write this term
$$ca_j^qJ_i^{(i+r-1)q} = cb_j^qJ_i^{(i+r-1)q} - ct_j^qJ_i^{(i+r-1)q}.$$
Using (3.4) (note $i+r-1\leq w+\l$), we obtain
$$ca_j^qJ_i^{(i+r-1)q} \inc cb_j^qJ_i^{(i+r-1)q} + J_{j-1}^{(i+r-1)q}$$
and so
$$\multline
c^{i-g} J_i^{(i+r)q}\inc \\
c^{i-g-1} [cJ_g\brq J_i^{(i+r-1)q}
+ (c b_{g+1}^q J_i^{(i+r-1)q} + J_g^{(i+r-1)q})
+\cdots + (c b_{i}^q J_i^{(i+r-1)q} + J_{i-1}^{(i+r-1)q})],
\endmultline$$
which by the induction hypothesis is contained in
$$
J_g\brq (B_i^r)\brq + b_{g+1}^q (B_i^r)\brq + (B_b^{r+i-g})\brq
+ \cdots + b_i^q (B_i^r)\brq + (B_{i-1}^{r+1})\brq \inc (B_i^{r+1})\brq.
$$

In particular, note that $$c^{\l-g} J_\l^{(\l+r)q} \inc  (B_\l^{r+1})\brq \tag{3.5}$$
for all $r\leq w$.

We now prove that  $\ic {I^{\l+w}}\inc (\A^{w+1})^*$.
 Let $u\in \ic {I^{\l+w}}$. Choose an element  $d \in R^0$ such that
$du^q\inc J^{(\l+w) q}$. Then
$c^{\l-g}du^q\in c^{\l-g}J^{(\l+w) q}\inc (B_\l^{w+1})\brq$ by (3.5).
It follows that $u\in (B_\l^{w+1})^* = (\A^{w+1})^*$. \qed\enddemo

\remark{Remark}  Theorem 3.3 is still valid even if $ht(I) = 0$.  In
this case choose $c_1 \in I^M$ and $c_2$ in the intersection
of the minimal primes of $0$ which do not contain $I$ and avoiding those that contain I.
Thus $c_2 I^N = 0$ for $N \gg 0$ and $c = c_1 + c_2 \in R^0$
satisfies equation (3.4).
\endremark

An almost immediate consequence is one of our main theorems:

\proclaim{Theorem 3.6}  Let $(R,\m)$ be an $F$-rational local ring of
positive characteristic $p$,
and let $I \inc R$ be an ideal generated by $\l$ elements.
Then $\ic{I^{\l+w}}
\inc I^{w+1}$ for all $w \ge 0$.
\endproclaim
 
\demo{Proof}
There is no loss of generality in assuming that $R$ has an
infinite residue field. We can replace $I$ by a minimal reduction of
itself;  suppose that $J$ 
is that minimal reduction. The number of generators of $J$ is at most
$\ell$, so without loss of generality we may assume $\ell$ is
the number of generators of $J$. Fix an integer $N$. We think of $w$ as fixed, and
choose $t_{g+1},\ldots, t_\l$ and $a_1,...,a_{\l}$ as in Proposition 3.2. 
In particular, $t_i\in \m^N$ for all $i$.
By (3.3), $\ic{I^{\l+w}} \inc (\A^{w+1})^* = \A^{w+1} 
\inc J^{w+1}+(t_{h+1},\ldots, t_\l)\inc J^{w+1} + \m^N$.
The equality $(\A^{w+1})^* = \A^{w+1}$ above follows from \cite{A, Thm. 1.1}.
By the Krull intersection theorem we obtain that
$\ic{I^{\l+w}} \inc \cap_N (J^{w+1}+\m^N) = J^{w+1}$.  \qed
\enddemo

This characteristic $p$ theorem allows us to prove the same result in
equicharacteristic $0$:

\proclaim{Theorem 3.7} Let $R$ be an algebra of finite type over a
field of characteristic $0$ and having only rational singularities.
Let  $I \inc R$ be an ideal generated by  $\l$ elements.
Then $\ic{I^{\l+w}}
\inc I^{w+1}$ for all $w\geq 0$.
\endproclaim

\demo{Proof}
By the work of Hara \cite{Ha} and independently Mehta and
Srinivas \cite{MS}, $R$ is of F-rational type. It is straightforward to prove 
in this case that if the conclusion holds in a dense open set of fibers
in some family of models $A\ra R_A$ of $R$, it also holds in $R$. Hence
we may pass to positive characteristic and assume that $R$ is finitely
generated over a field of characteristic $p > 0$ such that $R_P$ is
F-rational for all primes $P$. The conclusion will follow if we prove
it locally as the number of generators can only drop after localization.
It follows that we can reduce to the local F-rational case, and apply
Theorem 3.6 to finish the proof.  \qed

\enddemo
\bigskip
\head{4. F-Rational Gorenstein Rings}
\endhead
\bigskip

Our next theorem is new, even in the case $R$ is regular, as far as
we know. The proof is based on a careful analysis of the proof of
Theorem 3.5, and  the ideas behind the cancellation
theorem of \cite{Hu}. See also \cite{CP} for further cancellation
results. Our main theorem applies to rings which are F-rational and
Gorenstein. It is known \cite{HH3, (3.4), (4.7)}  that F-rational and F-regular are the same
when the base ring is Gorenstein. A ring $R$ is \it F-regular \rm if
every ideal is tightly closed in every localization of $R$. Of course,
all regular rings are F-regular, but 
the class of F-regular rings is considerably broader than that of regular rings.

\proclaim{Theorem 4.1}  Let $(R,\m)$ be an F-rational Gorenstein local ring
of dimension $d$ and having positive characteristic.
Suppose that $I$ is an ideal of height $g$, analytic spread $\l > g$
with $R/I$ CM.  For any reduction $J$ of $I$, $\ic{I^{\l-1}} \inc J$.
\endproclaim
 
\demo{Proof}
There is no loss of generality in assuming that $R$ has an infinite
residue field and that $J$ is a minimal reduction.
Fix an integer $N$ and set $w = 0$ in the notation of Proposition 3.2 and Theorem
3.3.
 We will prove that $\ic{I^{\l-1}} \inc J + \m^N$. An
application of the Krull intersection theorem then finishes the proof.

 We choose $t_{g+1},\ldots, t_\l$
and $a_1,...,a_{\l}$
as in Proposition 3.2, with $N$ fixed as above.
Let $b_i = a_i + t_i$ for $1 \le i \le \l$.
Choose $\bx = x_{\l+1},\ldots, x_d$ so that $b_{g+1}, \ldots, b_\l,
\bx$ is a regular sequence on $R/I$ and set $\A = (\vec b \l, \bx)$.
We set $D = J_g:t_{g+1}$ and $K = (J_g,b_{g+2},\ldots, b_\l,\bx)$.

Let $Q = (I,b_{g+2},\ldots, b_\l, \bx) + K:D$.
We claim that $\A:t_{g+1} \inc Q$.  Suppose that
$$t_{g+1} u = w + v b_{g+1} \tag{4.2}$$ where $w \in K$.
Then $t_{g+1}(u-v) \in (J_{g+1},b_{g+2},\ldots, b_\l,\bx)$ and
hence
$$
u -v \in (J_{g+1},b_{g+2},\ldots, b_\l,\bx):b_{g+1} \inc
(I,b_{g+2},\ldots, b_\l,\bx):b_{g+1} \inc (I,b_{g+2},\ldots, b_\l,\bx)
$$
since $R/I$ is Cohen-Macaulay. Hence $u-v\in Q$ and to prove $u\in Q$
it suffices to prove that $v\in K:D$. Let $d\in D$ and consider $dv$.
Using (4.2) we obtain that $t_{g+1}du = dw + dv b_{g+1}$, and hence
$dvb_{g+1}\in   
(J_g,b_{g+2},\ldots, b_\l,\bx)$. Thus  $$D_ v  \inc
(J_g,b_{g+2},\ldots, b_\l,\bx):b_{g+1} = (J_g,b_{g+2},\ldots, b_\l,\bx) = K.$$
This proves our claim, and in particular proves that $\A:Q \inc \A:(\A:t_{g+1})$.
 
We next claim that $\ic{I^{\l-1}} \inc \A:Q$. First observe that
$(I,b_{g+2},\ldots, b_\l, \bx)\cdot\ic{I^{\l-1}}\inc I\cdot\ic{I^{\l-1}} + \A$,
and by Theorem 2.6, $I\cdot\ic{I^{\l-1}}\inc \A$ (using that $R$ is F-rational).
Hence it remains to prove that $\ic{I^{\l-1}}\cdot (K:D)\inc \A$. We use a lemma.

\proclaim{Lemma 4.3} With the notation as above,
$$t_{g+1}\cdot\ic{I^{\l-1}}\inc J_{g}.$$
\endproclaim

\demo{Proof of Lemma 4.3} Let $z\in \ic{I^{\l-1}}$ and choose an element
$d\in R^o$ so that $dz^n\in I^{n(\l-1)}$ for all $n$. Choose $c\in I^M$ nonzero
as in (3.4). Using (3.4) we then obtain that
$$dct_{g+1}^qz^q\in ct_{g+1}^qI^{q(\l-1)}\inc t_{g+1}^qI^{q(\l-1)+M}\inc J_g^{q(\l-1)}\inc
J_g\brq$$
the last containment following as $\l-1\geq g$ and $J_g$ has $g$ generators.
Hence $t_{g+1}z\in (J_g)^* $. As $R$ is F-rational  $t_{g+1}z\in J_g$, proving the lemma.\qed
\enddemo

The Lemma proves that $\ic{I^{\l-1}} \inc
D$.  Hence $\ic{I^{\l-1}}((J_g,b_{g+2},\ldots, b_\l,\bx):D) \inc \A$.
We have proved that $\ic{I^{\l-1}} \inc \A:Q$.
 
By local duality, we have $\ic{I^{\l-1}} \inc \A:Q\inc  \A:(\A:t_{g+1})
\inc (J_{g+1}, t_{g+1},b_{g+2},\ldots, b_\l,\bx) \inc
(J,t_{g+1},\ldots, t_\l,\bx)\inc J + \m^N$. \qed 
\enddemo

\demo{Acknowledgement} The authors thank Reinhold H\"ubl for valuable 
corrections in our original preprint.
\enddemo

\enddocument